\documentclass[12pt]{article} 
\usepackage[utf8]{inputenc}
\usepackage[russian,english]{babel}
\usepackage{amsfonts}
\usepackage{amstext}
\usepackage{graphicx}
\usepackage{amsmath}
\usepackage{amssymb}
\usepackage{euscript}
\usepackage{amsthm}
  
\graphicspath{{Graphics/}}
\newtheorem{lemma}{Lemma}
\newtheorem{theorem}{Theorem}
\newtheorem{prop}{Proposition}
\newtheorem{rem}{Remark}

\begin{document}

\title{Initial value problem for the
 time-dependent linear Schr\"odinger equation 
with a point singular potential
 by the unified transform method}
\author{Ya. Rybalko \\
 \small\em V.N.Karazin Kharkiv National University\\
 \small\em B.Verkin Institute for Low Temperature Physics and Engineering}

\date{}

\maketitle

\begin{abstract}
We study an initial value problem for the one-dimensional non-stationary linear Schr\"odinger equation with a point singular potential. 
In our approach,  the problem 
is considered as a  system 
of coupled initial-boundary value (IBV) problems on two half-lines,
to which we apply the unified approach 
to IBV problems for linear and integrable nonlinear equations,
also known as the Fokas unified transform 
method.
Following the ideas of this method, we obtain the integral 
representation of the solution of the initial value problem.
\end{abstract}

\section{Introduction}
\label{intr}

The nonlinear Schr\"odinger (NLS) equation with an external potential $V(x)$ 
(the so-called Gross-Pitaevskii (GP) equation)
\begin{equation}
\label{1}
iu_{t}+V(x)u+\bigtriangleup u+2\kappa u|u|^{2}=0,\qquad \kappa=\pm1
\end{equation}
is used to describe a lot of phenomena in physics. In particular, it describes the static and dynamical properties of Bose-Einstein condensates, which attracts considerable interest of the researchers after it was experimentally observed in 1995 (see e.g. \cite{KFC} and the references therein).
Also the GP equation is widely used for  modelling of superconductors, for describing
optical vortices that resemble small twisters in a superfluid, and it appears in the studies of laser beams in Kerr media and focusing nonlinearity \cite{R}.
For applications it is desirable to get exact solution of the initial value (IV) problem for equation (\ref{1}), because it is helpful for describing detailed aspects of a particular physical system, when the approximate methods could be inadequate. But even in one dimension, it is hard to obtain such solutions of the original nonlinear problem. 

A reduction of (\ref{1}) related to short-range interactions
consists in replacing $V(x)$ by a point singular potential, which, in the one-dimensional case, is given by $V(x)=q\delta(x)$, $q\in\mathbb{R}$.
Assuming that the initial data 
$u(x,0)=u_0(x)$ is an {\em even} function, the solution is even for all $t$, 
and the IV problem

\begin{equation}
\label{1.0}
\left\{
\begin{array}{lcl}
iu_{t}+u_{xx}+q\delta(x)u+2u|u|^{2}=0,& x>0, \ t>0,\\
u(x,0)=u_{0}(x),& x\ge 0.\\
\end{array}
\right. 
\end{equation}
can be reduced 
 to the initial-boundary value (IBV) problem
for the NLS equation without the singular term, with the Robin {\em homogeneous}
boundary condition, see \cite{DP}:
\begin{equation}
\label{1.1}
\left\{
\begin{array}{lcl}
iu_{t}+u_{xx}+2u|u|^{2}=0,& x>0, \ t>0,\\
u_{x}(0,t)+qu(0,t)=0,& t\ge 0,\\
u(x,0)=u_{0}(x),& x\ge 0.\\
\end{array}
\right. 
\end{equation}

In the latter problem, the  boundary conditions are called \textit{linearizable},
because in this case there exists an adaptation of the Inverse Scattering Transform (IST) method, which turns to be as efficient as the IST method for the 
initial value problem (on the whole line) (see e.g. \cite{FIS,SI}). In other words, an initial boundary value problem with linearizable boundary conditions
is integrable: the appropriate version of the IST method 
 reduces it  to a series of {\em linear} problems in a way similar to the problems on the whole line, which in turn allows obtaining 
detailed information about properties of the solution, e.g., its the long-time behavior.

Particularly, in \cite{DP} (see also \cite{HZ}) the authors have studied in details 
the long-time behavior of the solution of (\ref{1.1}) by applying the nonlinear steepest-descent method
for Riemann--Hilbert problems (originally introduced for the whole line problem \cite{DIZ,DZ}).

However, if $u_0(x)$ is not assumed to be even,
(\ref{1.0}) does not reduce to (\ref{1.1}) and thus the approach
of \cite{DP} cannot be applied directly.

When analyzing a nonlinear problem, it is natural to solve, first, the associated linearized problem:
\begin{equation}
\label{1.2}
\left\{
\begin{array}{lcl}
iu_{t}+q\delta(x-a)u+u_{xx}=0,& x\in\mathbb{R}, \ t>0,\\
u(x,0)=u_{0}(x), & x\in\mathbb{R}. \\
\end{array}
\right. 
\end{equation}
Problem (\ref{1.2}) is the main object of study in the  present paper.
When trying to reduce (\ref{1.2}) to a problem on a half-line, we arrive, instead of  (\ref{1.1}),
at a system of \textit{two coupled half-line problems}, see (\ref{d3}) and (\ref{d4}) below. In other words, we treat the problem (\ref{1.2}) as a so-called \textit{interface problem}, with one interface between two infinite domains. Recall that the interface problems are initial-boundary value problems, for which the boundary conditions depend on the solution of the problems in adjacent domains \cite{DPS}. These problems have  recently attracted considerable attention,
see, e.g., \cite{DS, SD, DSS, DPS}, where the authors obtained explicit solutions of the original problems by applying the approach 
known as the Fokas unified transform method (initiated by A.Fokas \cite{F} and successfully developed in the further works, see, e.g.,  \cite{DTV, F1, FP}). 

The unified transform method is based on the so-called \textit{Lax pair} representation of an equation in question. Recall that the Lax pair of a given (partial differential) equation is a system of linear ordinary differential equations involving an additional (spectral) parameter, such that the given equation is the compatibly condition of this system. Originally, the Lax pair representation was introduced for certain nonlinear evolution equations, called integrable (see e.g. \cite{AKNS}). This representation plays a key role in solving \textit{ initial value} problems for such equations by the inverse scattering transform (IST) method, where one (spatial) equation from the pair establishes a change of variables, passing from functions of the spatial variable to functions of the spectral parameter, whereas the evolution in time turns out to be, due to the other (temporal) equation from the Lax pair, linear. The Lax pair representation also plays a crucial role in studying {\em initial  boundary value} problems, where, according to the  unified transform method, both equations from the Lax pair are treated in a same manner,  as spectral problems.

While only very special (although very important) nonlinear equations have a Lax pair representation (the NLS equation in (\ref{1.1}) is one of them), the Lax pair for linear equations with constant coefficients can be constructed algorithmically. Particularly, the Lax pair for the linearized form of the NLS equation,
 $iu_{t}+u_{xx}=0$, is as follows:
\begin{equation}
\label{2.5}
\left\{
\begin{array}{lcl}
(\mu e^{-ikx+ik^{2}t})_{x}=ue^{-ikx+ik^{2}t}\\
(\mu e^{-ikx+ik^{2}t})_{t}=(iu_{x}-ku)e^{-ikx+ik^{2}t}\\
\end{array}
\right.
\end{equation}
Indeed, 
direct calculations show that the compatibly condition of (\ref{2.5}),
which is the equality 
$$(ue^{-ikx+ik^{2}t})_{t}=((iu_{x}-ku)e^{-ikx+ik^{2}t})_{x}
$$
 to be satisfied for all $k\in{\mathbb C}$, reads  $iu_{t}+u_{xx}=0$. 

In this paper we obtain an integral representation of the solution of problem (\ref{1.2}) using the ideas of the Fokas method (see Theorem \ref{dutv}).
Since this method for both linear and integrable nonlinear equations uses
the similar ingredients --- the Lax pair representation and the analysis
of the so-called global relation (see below) --- we believe that our results can be useful
for studying the initial value problem for the  {\em nonlinear} Schr\"odinger equation with a point singular potential.

\section{Representation of the solution}
\label{oned}

In what follows we will consider 
 the initial value problem (\ref{1.2}),
where $u_0(x)$ is a smooth function decaying as $|x|\to \infty$, $a\in\mathbb{R}$ is a parameter,
and the solution $u(x,t)$ is assumed to decay to $0$ as $|x|\to \infty$ for all $t>0$ .

First, let us assume that the solution $u(x,t)$ of the problem (\ref{1.2}) exists,
such that $u(x,t)\in L^{1}({\mathbb R})$ and $u_{x}(x,t)\rightarrow 0$ as $|x|\rightarrow\infty$ for all $t$. 
Our goal is to obtain the (integral) representation for $u(x,t)$ in terms of the given data $u_0(x)$.

In order to give exact meaning for (\ref{1.2}), we notice that the delta function $2q\delta(x-a)$ is to be understood as 
introducing the jump condition for the first $x$-derivative:
\begin{equation}\label{u-jump-deriv}
u_{x}(a+0,t)-u_{x}(a-0,t)+2qu(a,t)=0,\qquad t>0
\end{equation}
whereas $u(x,t)$ is continuous across $x=a$: 
\begin{equation}\label{u-contin}
u(a+0,t)=u(a-0,t), \qquad t>0.
\end{equation}
Thus
(\ref{1.2}) is to be understood as the following system of equations:
\begin{equation}
\label{d2}
\left\{
\begin{array}{lcl}
iu_{t}+u_{xx}=0, & x\in(-\infty,a)\cup (a,\infty), t>0,\\
u(a+0,t)=u(a-0,t),& t>0 \\
u_{x}(a+0,t)-u_{x}(a-0,t)+2qu(a,t)=0,& t>0,\\
u(x,0)=u_{0}(x),& x\in(-\infty,\infty)\\
\end{array}
\right.
\end{equation}

Let's introduce the notations:
\begin{equation}
\label{d2.1}
\begin{matrix}
u^{(1)}(x,t):=
\left\{
\begin{array}{lcl}
u(x,t),& x\leq a\\
0,& x>a\\
\end{array}
\right.,
&
u^{(2)}(x,t):=
\left\{
\begin{array}{lcl}
u(x,t), & x\geq a\\
0, & x<a\\
\end{array}
\right.,
\end{matrix}
\end{equation}

\begin{equation}
\label{d2.2}
\begin{matrix}
u^{(1)}_{0}(x):=
\left\{
\begin{array}{lcl}
u_{0}(x),& x\leq a\\
0,& x>a\\
\end{array}
\right.,
&
u^{(2)}_{0}(x):=
\left\{
\begin{array}{lcl}
u_{0}(x),& x\geq a\\
0,& x<a\\
\end{array}
\right..\\
\end{matrix}
\end{equation}
Then the jump relation can be written as
$$u^{(1)}_{x}(a,t)-qu^{(1)}(a,t)=u^{(2)}_{x}(a,t)+qu^{(2)}(a,t)=:\phi(t),
$$
which introduces $\phi(t)$ for $t>0$.
On the other hand, the continuity of $u(x,t)$ allows introducing 
$$g(t):=u^{(1)}(a,t)=u^{(2)}(a,t).$$

Therefore,  the initial value  problem (\ref{d2}) can be  equivalently written as a \textit{ pair of coupled initial boundary value problems}, 
with the   inhomogeneous  boundary conditions:
\begin{subequations}
\begin{equation}
\left\{
\begin{array}{lcl}
\label{d3}
iu_{t}^{(1)}+u_{xx}^{(1)}=0,& x<a,\,t>0,\\
u^{(1)}(x,0)=u_{0}^{(1)}(x),& x<a,\\
u_{x}^{(1)}(a,t)-qu^{(1)}(a,t)=\phi(t),& t>0,\\
u^{(1)}(a,t)=g(t),& t>0.\\
\end{array}
\right.
\end{equation}
\begin{equation}
\left\{
\begin{array}{lcl}
\label{d4}
iu_{t}^{(2)}+u_{xx}^{(2)}=0,& x>a,\,t>0,\\
u^{(2)}(x,0)=u_{0}^{(2)}(x),& x>a,\\
u_{x}^{(2)}(a,t)+qu^{(2)}(a,t)=\phi(t),& t>0,\\
u^{(2)}(a,t)=g(t),& t>0.\\
\end{array}
\right.
\end{equation}
\end{subequations}
We emphasize  that neither $\phi(t)$ nor $g(t)$ are given as data for the problems,
but (\ref{d3}) and (\ref{d4}) are coupled by the condition that these
functions are the same for the both problems.

 In what follows, we 
will use the following notations for the direct and inverse 
Fourier transforms: 
$$\EuScript{F}(f)(k)=\hat{f}(k)=\int_{-\infty}^{\infty}e^{-ikx}f(x)\,dx,\,\EuScript{F}^{-1}(\hat{f})(x)=f(x)=\frac{1}{2\pi}\int_{-\infty}^{\infty}e^{ikx}\hat{f}(x)\,dx.$$
With this notations, we define
\begin{equation}\label{u-hat}
\hat{u}^{(1)}_{0}(k):=\EuScript{F}(u^{(1)}_{0})(x) = \int_{-\infty}^{a}e^{-ikx}u_0(x)dx, \qquad 
\hat{u}^{(2)}_{0}(k):=\EuScript{F}(u^{(2)}_{0})(x) = \int_{a}^{\infty}e^{-ikx}u_0(x)dx.
\end{equation}

Now, let us assume for a moment that the functions $g(t)$ and $\phi(t)$
are given. Then, using the Lax pair representation separately for the both equations $iu^{(i)}_{t}+u^{(i)}_{xx}=0$, $i=1,2$, 
we can obtain explicitly
 the solutions of  (\ref{d3}) and (\ref{d4}). 

Introduce the Fourier-type transforms for the ``boundary values'' $g(t)$ and $\phi(t)$:
\begin{equation}
\label{d5}
\begin{array}{lcl}
h_{1}(k,t):=\int_{0}^{t}e^{ik^2(\tau-t)}g(\tau)\,d\tau,\qquad
h_{2}(k,t):=\int_{0}^{t}e^{ik^2(\tau-t)}\phi(\tau)\,d\tau.
\end{array}
\end{equation}

\begin{prop}
The solutions $u^{(j)}(x,t)$, $j=1,2$ of problems  (\ref{d3}) and (\ref{d4})
 can be obtained in terms of $h_{j}(k,t)$ and $\hat{u}^{(j)}_{0}(k)$, $j=1,2$
as follows:
\begin{equation}
\label{d6}
u^{(1)}(x,t)=\EuScript{F}^{-1}\left(e^{-ik^{2}t}\hat{u}^{(1)}_{0}(k)+e^{-ika}(iq-k)h_{1}(k,t)
+ie^{-ika}h_{2}(k,t)\right),
\end{equation}
\begin{equation}
\label{d7}
u^{(2)}(x,t)=\EuScript{F}^{-1}\left(e^{-ik^{2}t}\hat{u}^{(2)}_{0}(k)+e^{-ika}(iq+k)h_{1}(k,t)
-ie^{-ika}h_{2}(k,t)\right).
\end{equation}
\end{prop}
\begin{proof}
Consider the Lax pair for equation $iu^{(1)}_t+u^{(1)}_{xx}=0$:
\begin{equation}
\label{d7.5}
\left\{
\begin{array}{lcl}
(\mu e^{-ikx+ik^{2}t})_{x}=u^{(1)}e^{-ikx+ik^{2}t}\\
(\mu e^{-ikx+ik^{2}t})_{t}=(iu^{(1)}_{x}-ku^{(1)})e^{-ikx+ik^{2}t}.\\
\end{array}
\right.
\end{equation}
\\
According to  (\ref{d7.5}),  the 1-form 
$$
W_{1}(x,t,k):=(u^{(1)}e^{-ikx+ik^{2}t})\,dx+((iu_{x}^{(1)}-ku^{(1)})e^{-ikx+ik^{2}t})dt
$$
is exact:
\begin{equation}
W_{1}(x,t,k)=d\left[(\mu e^{-ikx+ik^{2}t})\right],\,\,x\leq a,t\geq 0.
\end{equation}
\\
It follows that $\oint_{\Gamma}W_{1}(y,\tau,k)=0$,
where 
$\Gamma$ is the rectangle with the  vertices $(-X,0)$,  $(a,0)$, $(a,t)$, and 
$(-X,t)$, with any $X>0$, and thus
\begin{eqnarray}
\label{d8}
0&=&\int_{\Gamma}W_{1}(y,\tau,k)=\oint_{\Gamma}(u^{(1)}e^{-ikx+ik^{2}\tau})\,dx+((iu_{x}^{(1)}-ku^{(1)})e^{-ikx+ik^{2}\tau})\,d\tau\\
\nonumber
&=&\int_{-X}^{a}e^{-ikx}u^{(1)}_{0}(x)\,dx+
\int_{0}^{t}e^{-ika+ik^{2}\tau}(iu^{(1)}_{x}(a,\tau)-ku^{(1)}(a,\tau))\,d\tau\\
\nonumber
&&-\int_{-X}^{a}e^{-ikx+ik^{2}t}u^{(1)}(x,t)\,dx
-\int_{0}^{t}e^{ikX+ik^{2}\tau}(iu^{(1)}_{x}(-X,\tau)-ku^{(1)}(-X,\tau))\,d\tau.
\end{eqnarray}
Letting $X\rightarrow+\infty$ in 
 (\ref{d8}) and taking into account the boundary conditions in (\ref{d3}) and that 
 the last term decays to 0 as $X\rightarrow+\infty$, equation (\ref{d8}) takes the form
\begin{equation}
\hat{u}^{(1)}_{0}(k)+\int_{0}^{t}e^{-ika+ik^{2}\tau}\left[i(\phi(\tau)+qg(\tau))-kg(\tau)\right]\,d\tau-
e^{ik^{2}t}\hat{u}^{(1)}(k,t)=0,
\end{equation}
which, due to the first term, is valid for all $k$ with $\Im k\ge 0$.
Using the notations (\ref{d5}), the last equation can be written as
\begin{equation}
\label{d9}
\hat{u}^{(1)}(k,t)=e^{-ik^{2}t}\hat{u}^{(1)}_{0}(k)+e^{-ika}(iq-k)h_{1}(k,t)+ie^{-ika}h_{2}(k,t),\quad \Im k\ge 0,
\end{equation}
which  implies (\ref{d6}).

The  IBV problem for ${u}^{(2)}(x,t)$ can  be solved in a  similar way, integrating the  form $W_2$, which is similar to $W_1$
with ${u}^{(1)}$ replaced by ${u}^{(2)}$, along the rectangle  with the  vertices $(X,0)$,  $(a,0)$, $(a,t)$, and 
$(X,t)$ and passing to the limit $X\rightarrow+\infty$. This gives
\begin{equation}
\hat{u}^{(2)}_{0}(k)-e^{ik^{2}t}\hat{u}^{(2)}(k,t)-
\int_{0}^{t}e^{-ika+ik^{2}\tau}\left[i(\phi(\tau)-qg(\tau))-kg(\tau)\right]\,d\tau=0, \quad \Im k\le 0.
\end{equation}
Using (\ref{d5}), the last equation can be written as
\begin{equation}
\label{d11}
\hat{u}^{(2)}(k,t)=e^{-ik^{2}t}\hat{u}^{(2)}_{0}(k)+e^{-ika}(iq+k)h_{1}(k,t)-ie^{-ika}h_{2}(k,t),\quad \Im k\le 0,
\end{equation}
which  implies (\ref{d7}).
\end{proof}

\begin{rem}
\label{r1}
In the framework of the Fokas method,
equations (\ref{d9}) and (\ref{d11}) are called 
the \emph{global relations}: they relate the ``boundary'' values of $u(x,\tau)$, taken at $x=a$, $\tau=0$, and $\tau=t$, in the spectral terms, i.e.,
 in the form of the corresponding Fourier-type transforms. The global relations involves given, for a well posed IBV problem, initial/boundary data 
(in our case, these are the initial data $u_0(x)$)
as well as unknown boundary values (in our case, these are $g(t)$ and $\phi(t)$). 
The central importance of the global relations is that they allow to \emph{characterize} the unknown
boundary values in terms of the  known (given) ones \cite{BFS03}.
 While for the IBV problems for nonlinear PDE with non-linearizable boundary conditions, this characterization
can be effectively used for studying various aspects of the problem \cite{BFS06, F2, FIS}
 but does not give a ``complete'' solution to the problem (the expression of the  solution in terms of the data for a well-posed problem), 
in the case of linear PDE this characterization is expected to help solving the 
IBV completely \cite{F1, FP}. Our goal in this paper is actually to demonstrate the latter for problem (\ref{d2})
(or, equivalently,  for the coupled problems (\ref{d3}) and (\ref{d4})).
\end{rem}

An important feature of the global relations  (\ref{d9}) and (\ref{d11}) is the properties of analyticity
and decay of their left-hand sides. For linear problems, the analyticity and decay follow directly from their
definitions as the Fourier transforms of functions supported on the associated
half-lines, which is  summarized in the following
\begin{lemma}
\label{dl1}
The functions $h_{i}(k,t)$, $i=1,2$, defined by (\ref{d5})
 are analytic in ${\mathbb C}$.
The function $\hat{u}^{(1)}(k,t)$ is analytic in ${\mathbb C}_{+}$ and $\hat{u}^{(2)}(k,t)$ is analytic in ${\mathbb C}_{-}$ for all $t\geq0$, where
\begin{equation}
\label{d12}
{\mathbb C}_{-}:=\{k\in{\mathbb C}:\Im k<0\}, \qquad{\mathbb C}_{+}=\{k\in{\mathbb C}:\Im k>0\}.
\end{equation}
  Moreover, $e^{-ika}\hat{u}^{(1)}(-k,t)$
and $e^{ika}\hat{u}^{(2)}(k,t)$ decay to 0 as $k\to\infty$ in 
$\bar{\mathbb C}_{-}$ (closure of ${\mathbb C}_{-}$).
\end{lemma}

Now observe that, by definition,  $h_{i}(k,t)=h_{i}(-k,t)$, $i=1,2$. This observation and Lemma \ref{dl1} allow rewriting (\ref{d9}) and (\ref{d11}) as a  
 system of  equations
\begin{equation}
\label{d13}
\left\{
\begin{array}{lcl}
\hat{u}^{(1)}(-k,t)-e^{-ik^{2}t}\hat{u}^{(1)}_{0}(-k)=e^{ika}(iq+k)h_{1}(k,t)+ie^{ika}h_{2}(k,t),\\
\hat{u}^{(2)}(k,t)-e^{-ik^{2}t}\hat{u}^{(2)}_{0}(k)=e^{-ika}(iq+k)h_{1}(k,t)-ie^{-ika}h_{2}(k,t),
\end{array}
\right.
\end{equation}
which holds for all $k\in\bar {\mathbb C}_{-}$.

Now suppose for a moment that the functions $\hat{u}^{(i)}(k,t)$, $i=1,2$
are known. Then (\ref{d13}) can be considered as a system of two 
algebraic equations for  two unknown functions $h_{i}(k,t)$. 
The solution of this system is given  for $k\in\bar {\mathbb C}_{-}$ by 
\begin{eqnarray}
\label{d14}
h_{1}(k,t)=\frac{e^{-ika}\left(\hat{u}^{(1)}(-k,t)-e^{-ik^{2}t}\hat{u}_{0}^{(1)}(-k)\right)
+e^{ika}\left(\hat{u}^{(2)}(k,t)-e^{-ik^{2}t}\hat{u}_{0}^{(2)}(k)\right)}{2(iq+k)},\\
\label{d15}
h_{2}(k,t)=\frac{e^{-ika}\left(\hat{u}^{(1)}(-k,t)-e^{-ik^{2}t}\hat{u}_{0}^{(1)}(-k)\right)
-e^{ika}\left(\hat{u}^{(2)}(k,t)-e^{-ik^{2}t}\hat{u}_{0}^{(2)}(k)\right)}{2i}.
\end{eqnarray}

Now we are at a position to formulate and prove the main representation result.
\begin{theorem}
\label{dutv}
The solution  $u(x,t)$ of problem (\ref{d2})  is given, in terms of the Fourier transforms of the initial data,
 as follows:
\begin{itemize}
\begin{subequations} \label{s}
\item if $q<0$, then
\begin{equation}\label{s1}
u(x,t)=\frac{1}{2\pi}\int_{-\infty}^{\infty}e^{ikx-ik^{2}t}
\left(\hat{u}_0^{(1)}(k)-
\frac{iqe^{-2ika}\hat{u}_0^{(1)}(-k)-k\hat{u}^{(2)}_{0}(k)}
{iq+k}\right)\,dk,
\end{equation}
for $x<a$, and 
\begin{equation}\label{s2}
u(x,t)=\frac{1}{2\pi}\int_{-\infty}^{\infty}e^{ikx-ik^{2}t}
\left(\hat{u}_0^{(2)}(k)-\frac{iqe^{-2ik a}\hat{u}_0^{(2)}(-k)+k\hat{u}_0^{(1)}(k)}{iq-k}\right)\,dk,
\end{equation}
for $x>a$.
\item if $q>0$, then
\begin{align}\label{s3}
\nonumber
u(x,t)=&\frac{1}{2\pi}\int_{-\infty}^{\infty}e^{ikx-ik^{2}t}
\left(\hat{u}_0^{(1)}(k)-
\frac{iqe^{-2ika}\hat{u}_0^{(1)}(-k)-k\hat{u}^{(2)}_{0}(k)}
{iq+k}\right)\,dk\\
&+qe^{q(x-a)+iq^{2}t}\left(\int_{-\infty}^{a}e^{q(y-a)}u_{0}(y)\,dy+
\int_{a}^{\infty}e^{-q(y-a)}u_{0}(y)\,dy\right),
\end{align}
for $x<a$, and 
\\
\begin{align}\label{s4}
\nonumber
u(x,t)=&\frac{1}{2\pi}\int_{-\infty}^{\infty}e^{ikx-ik^{2}t}
\left(\hat{u}_0^{(2)}(k)-\frac{iqe^{-2ik a}\hat{u}_0^{(2)}(-k)+k\hat{u}_0^{(1)}(k)}{iq-k}\right)\,dk\\
& +qe^{-q(x-a)+iq^{2}t}\left(\int_{-\infty}^{a}e^{q(y-a)}u_{0}(y)\,dy+
\int_{a}^{\infty}e^{-q(y-a)}u_{0}(y)\,dy\right),
\end{align}
for $x>a$.
\end{subequations}
\end{itemize}
\end{theorem}

\begin{rem}
The difference in the form of the solution in the cases $q>0$ and $q<0$ is related to the fact that the denominator in (\ref{d14}) 
 has a zero (at $k=-iq$) in the domain $ {\mathbb C}_{-}$ of validity of (\ref{d14})  in the case $q>0$ only.
\end{rem}

\begin{proof}

Since $h_{i}(k,t)=h_{i}(-k,t)$,   the first equation in (\ref{d13})
gives the system
\begin{equation}
\label{d17}
\left\{
\begin{array}{lcl}
e^{-ika}\left(\hat{u}^{(1)}(-k,t)-e^{-ik^{2}t}\hat{u}^{(1)}_{0}(-k)\right)=(iq+k)h_{1}(k,t)+ih_{2}(k,t),\\
e^{ika}\left(\hat{u}^{(1)}(k,t)-e^{-ik^{2}t}\hat{u}^{(1)}_{0}(k)\right)=(iq-k)h_{1}(k,t)+ih_{2}(k,t),\\
\end{array}
\right.
\end{equation}
that holds for $k\in\mathbb R$, which,   getting rid of $h_2(k,t)$, leads to the equation
\begin{equation}
\label{d18}
e^{-ika}\hat{u}^{(1)}(-k,t)-e^{ika}\hat{u}^{(1)}(k,t) + e^{ika-ik^{2}t}\hat{u}^{(1)}_{0}(k)-
e^{-ika-ik^{2}t}\hat{u}^{(1)}_{0}(-k)=2kh_{1}(k,t),\quad k\in\mathbb R.
\end{equation}

Multiply (\ref{d18}) by $\frac{1}{2\pi}e^{ik(x-a)}$, integrate over $k$ from $-N$ to $N$, and pass to the limit $N\to\infty$. Then for $x<a$, the first terms
vanishes due to Lemma \ref{dl1} and Jordan's lemma (related to a half-circle in the lower half-plane of $k$)
 whereas the second term becomes $-u^{(1)}(x,t)$ (as the result of the  inverse Fourier transform),
and thus the resulting equation is
\begin{eqnarray}
\label{d20}
\nonumber
&&-u^{(1)}(x,t)+\frac{1}{2\pi}\int_{-\infty}^{\infty}e^{ikx-ik^{2}t}\hat{u}^{(1)}_{0}(k)\,dk-
\frac{1}{2\pi}\int_{-\infty}^{\infty}e^{ik(x-2a)-ik^{2}t}\hat{u}^{(1)}_{0}(-k)\,dk\\
&&=\frac{1}{2\pi}\int_{-\infty}^{\infty}e^{ik(x-a)}2kh_{1}(k,t)\,dk,\,\,x<a.
\end{eqnarray}

Now, in order to obtain $u^{(1)}(x,t)$ in terms of the given (initial) data, we have to express 
$I_1(x,t):=\frac{1}{2\pi}\int_{-\infty}^{\infty}e^{ik(x-a)}2kh_{1}(k,t)\,dk$ for $x<a$  in terms of $u_{0}^{(i)}(x)$. 
First, we notice that  $h_{1}(k,t)$ decays to 0 
as $k\to\infty$ in the quadrant $\Re k>0$, $\Im k <0$. Thus, by
Jordan's lemma related to a part of large circle in this quadrant,
the integral over the real axis can be deformed to that over  the contour $\gamma$, which is shown in Fig. \ref{dd1}:
\begin{equation}\label{i1}
I_1(x,t)=
\frac{1}{2\pi}\int_{\gamma}e^{ik(x-a)}2kh_{1}(k,t)\,dk,
\end{equation}
Substituting  (\ref{d14}) into the r.h.s. of (\ref{i1}),
we have
$$
I_1(x,t)=I_2(x,t)-I_3(x,t),
$$
where
$$
I_2(x,t) = \frac{1}{2\pi}\int_{\gamma}e^{ik(x-a)}k\frac{e^{-ika}
\hat{u}^{(1)}(-k,t)+
e^{ika}\hat{u}^{(2)}(k,t)}{iq+k}\,dk,
$$
and
\begin{equation}\label{i3}
I_3(x,t) =\frac{1}{2\pi}\int_{\gamma}e^{ik(x-a)-ik^{2}t}k\frac{e^{-ika}\hat{u}_{0}^{(1)}(-k)+
e^{ika}\hat{u}_{0}^{(2)}(k)}{iq+k}\,dk.
\end{equation}
Then Lemma \ref{dl1} and Jordan's lemma related to a part of big circle in the quadrant  $\Re k<0$, $\Im k <0$
 implies that $I_2(x,t)=0$.

The integral in the r.h.s. of (\ref{i3}) already gives $I_3$ (and thus $I_1$) in terms of the initial  data only.
But now we can deform back the integration pass, from $\gamma$ to the real axis,
again using Lemma \ref{dl1} and Jordan's lemma related to a  part of big circle in the quadrant  $\Re k>0$, $\Im k <0$,
which gives:
\begin{itemize}
\item
if $q<0$, then 
$$
I_{3}(x,t) = \frac{1}{2\pi}\int_{-\infty}^{\infty}e^{ik(x-a)-ik^{2}t}k\frac{e^{-ika}\hat{u}_{0}^{(1)}(-k)+e^{ika}\hat{u}_{0}^{(2)}(k)}{iq+k}\,dk.
$$
\item
if $q>0$, then 
$$
I_3(x,t) = I_{31}(x,t)+I_{32}(x,t),
$$
where 
\begin{align*}
I_{31}(x,t) & = i\underset{k=-iq}{\operatorname{Res}}e^{ik(x-a)-ik^{2}t}k\frac{e^{-ika}\hat{u}_{0}^{(1)}(-k)+
e^{ika}\hat{u}_{0}^{(2)}(k)}{iq+k} \\
& = qe^{q(x-a)+iq^{2}t}\left(\int_{-\infty}^{a}e^{q(y-a)}u_{0}^{(1)}(y)\,dy+
\int_{a}^{\infty}e^{-q(y-a)}u_{0}^{(2)}(y)\,dy\right),
\end{align*}
and
$$
I_{32}(x,t) = \frac{1}{2\pi}\int_{-\infty}^{\infty}e^{ik(x-a)-ik^{2}t}k\frac{e^{-ika}\hat{u}_{0}^{(1)}(-k)+e^{ika}\hat{u}_{0}^{(2)}(k)}{iq+k}\,dk.
$$
\end{itemize}
Substituting this into (\ref{d20})
 we arrive at the statements of Theorem \ref{dutv}
concerning $u(x,t)$ for $x<a$.

\begin{figure}[h]
\begin{minipage}[h]{0.5\linewidth}
\centering{\includegraphics[width=0.36\linewidth]{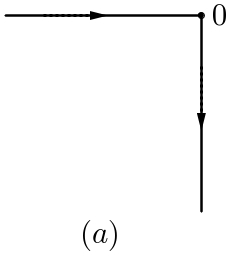}}
\end{minipage}
\begin{minipage}[h]{0.5\linewidth}
\centering{\includegraphics[width=0.4\linewidth]{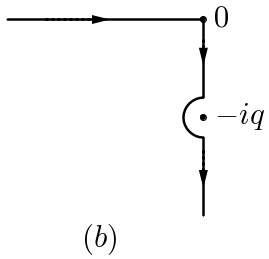}}
\end{minipage}
\caption{Contour $\gamma$: $(a)$ in the case $q<0$; $(b)$ in the case $q>0$.}
\label{dd1}
\end{figure}

To find $u^{(2)}(x,t)$, we proceed in the similar way, 
with the starting point being the equation (cf. (\ref{d18})), which follows from  
 the second equation in (\ref{d13}):
\begin{equation}
\label{d23}
e^{ika}\hat{u}^{(2)}(k,t)-e^{-ika}\hat{u}^{(2)}(-k,t)+
e^{-ika-ik^{2}t}\hat{u}^{(2)}_{0}(-k)-e^{ika-ik^{2}t}\hat{u}^{(2)}_{0}(k)=2kh_{1}(k,t), \quad k\in\mathbb R.
\end{equation}
Multiply (\ref{d23}) by $\frac{1}{2\pi}e^{ikx}$, integrate over $k$ from $-N$ to $N$, and pass to the limit $N\to\infty$. Then for $x<0 $, the first terms
vanishes due to Lemma \ref{dl1} and Jordan's lemma (related to a half-circle in the lower half-plane of $k$)
 whereas the second term becomes  $-u^{(2)}(a-x,t)$ (again as the result of the inverse Fourier transform):
\begin{eqnarray*}
\label{d24}
&&-u^{(2)}(a-x,t)+\frac{1}{2\pi}\int_{-\infty}^{\infty}e^{ik(x-a)-ik^{2}t}\hat{u}^{(2)}_{0}(-k)\,dk-
\frac{1}{2\pi}\int_{-\infty}^{\infty}e^{ik(x+a)-ik^{2}t}\hat{u}^{(2)}_{0}(k)\,dk\\
&&=\frac{1}{2\pi}\int_{-\infty}^{\infty}e^{ikx}2kh_{1}(k,t)\,dk,\quad x<0,
\end{eqnarray*}
or
\begin{eqnarray}
\label{d25}
\nonumber
&&-u^{(2)}(x,t)+\frac{1}{2\pi}\int_{-\infty}^{\infty}e^{-ikx-ik^{2}t}\hat{u}^{(2)}_{0}(-k)\,dk-
\frac{1}{2\pi}\int_{-\infty}^{\infty}e^{-ik(x-2a)-ik^{2}t}\hat{u}^{(2)}_{0}(k)\,dk\\
&&=\frac{1}{2\pi}\int_{-\infty}^{\infty}e^{-ik(x-a)}2kh_{1}(k,t)\,dk,\quad x>a.
\end{eqnarray}
Taking into account that the last term in (\ref{d25}) has already been 
calculated above, from (\ref{d25}) we obtain
\begin{equation}
\label{d26}
u^{(2)}(x,t)=\frac{1}{2\pi}\int_{-\infty}^{\infty}e^{-ikx-ik^{2}t}\hat{u}^{(2)}_{0}(-k)\,dk-
\frac{1}{2\pi}\int_{-\infty}^{\infty}e^{-ikx-ik^{2}t}\frac{iqe^{2ika}\hat{u}^{(2)}_{0}(k)-
k\hat{u}^{(1)}_{0}(-k)}{iq+k}\,dk,
\end{equation}
in the case $q<0$ and 
\begin{eqnarray}
\nonumber
\label{d27}
u^{(2)}(x,t)&=&\frac{1}{2\pi}\int_{-\infty}^{\infty}e^{-ikx-ik^{2}t}\hat{u}^{(2)}_{0}(-k)\,dk-
\frac{1}{2\pi}\int_{-\infty}^{\infty}e^{-ikx-ik^{2}t}\frac{iqe^{2ika}\hat{u}^{(2)}_{0}(k)-
k\hat{u}^{(1)}_{0}(-k)}{iq+k}\,dk\\
&&
+ qe^{-q(x-a)+iq^{2}t}\left(\int_{-\infty}^{a}e^{q(y-a)}u_{0}^{(1)}(y)\,dy+
\int_{a}^{\infty}e^{-q(y-a)}u_{0}^{(2)}(y)\,dy\right),
\end{eqnarray}
in the case $q>0$.
Substituting $k$ by $-k$ in (\ref{d26}) and (\ref{d27})
we arrive at the statements of the Theorem \ref{dutv}
concerning $u(x,t)$ for $x>a$.
\end{proof}

\begin{rem}
Setting $x=a$ in (\ref{s}) gives one part of the ``generalized Dirichlet-to-Neumann map'',
which has to give unknown boundary values (in the  case of problems  (\ref{d3}) and (\ref{d4}), they are  $g(t)$ and $\phi(t)$)
 in terms of the given data (in our case, they are the  initial data $u_0(x)$):
\begin{subequations}\label{b}
\begin{equation}
\label{b1}
g(t)=\frac{1}{2\pi}\int_{-\infty}^{\infty}e^{-ik^{2}t}
k\frac{e^{-ika}\hat{u}_0^{(1)}(-k)+e^{ika}\hat{u}^{(2)}_{0}(k)}
{iq+k}\,dk,
\end{equation}
for $q<0$, and
\begin{align}
\nonumber
g(t)=&\frac{1}{2\pi}\int_{-\infty}^{\infty}e^{-ik^{2}t}
k\frac{e^{-ika}\hat{u}_0^{(1)}(-k)+e^{ika}\hat{u}^{(2)}_{0}(k)}
{iq+k}\,dk\\
\label{b2}
&+qe^{iq^2t}
\left(\int_{-\infty}^{a}e^{q(y-a)}u_{0}(y)\,dy+
\int_{a}^{\infty}e^{-q(y-a)}u_{0}(y)\,dy\right),
\end{align}
for $q>0$.
\end{subequations}

 In order to obtain the second part of this map, i.e., $\phi(t)$  in terms of  $u_0(x)$,
one can  use equation (\ref{d15}) directly following from the global relations.
Indeed, multiply (\ref{d15}) by $\frac{k}{\pi}e^{-ik^2 t'}$ with $t'<t$, integrate along $\gamma$ and take $t'\to t$. Then the l.h.s. gives $\phi(t)$
whereas in the r.h.s, integrating by parts, interchanging the integrals and calculating explicitly the integral $\int_{\gamma} e^{-ikx-ik^2 t} dk$,
one arrives at the representation

\begin{equation}\label{b3}
\phi(t)=\frac{e^{-i\frac{\pi}{4}}}{2\sqrt{\pi t}}
\left(\int_{-\infty}^{a} + \int_{a}^{\infty}\right)
e^{\frac{i(x-a)^{2}}{4t}}u_0^{\prime}(x)\,dx,\quad q\in\mathbb{R}.
\end{equation}
for all $q\in\mathbb{R}$.

\end{rem}

\begin{rem}
\label{prop0}
Assuming that the initial data $u_0(x)$ is a smooth function for $x\in(-\infty,a)\cup(a,\infty)$ decaying fast as $|x|\to \infty$ 
and  matching well the ``inner'' conditions  (\ref{u-jump-deriv}) and  (\ref{u-contin}), i.e.,
\begin{equation}\label{u0-cond}
u_0(a-0)=u_0(a+0), \quad u'_0(a+0)- u'_0(a-0)+2qu_0(a)=0, \quad  u''_0(a+0)= u''_0(a-0)
\end{equation}
(the latter condition is suggested by  (\ref{u-contin}) and the equation $iu_t+u_{xx}=0$ considered at $x=a\pm 0$)
formulas (\ref{s}) in Theorem 1 give a classical solution to  (\ref{d2}).
\end{rem}
Indeed, assume that 
 $u_{0}(x)\in C^{4}({\mathbb R}\backslash\{a\})\cap C(\mathbb{R})$,
$\frac{d^p}{dx^p}u_{0}(x)\in L^{1}({\mathbb R}\setminus\{a\})$ ($p=0,\dots,4$),
 $x^j\frac{d^4}{dx^4}u_0(x)\in L^{1}(\mathbb{R}\setminus\{a\})$ ($j=1,2$), and $u_0(x)$ satisfies (\ref{u0-cond}).
Then, integrating by parts the definitions (\ref{u-hat}) of $\hat{u}_0^{(j)}(k)$ and taking into account (\ref{u0-cond}),
one finds that the  integrands in (\ref{s}) are $O(\frac{1}{k^4})$ as $k\to\infty$. 
It follows that  $u(x,t)$ defined by (\ref{s}) has  partial derivatives $u_t(x,t)$ and $u_{xx}(x,t)$
in the domains $x<a, t>0$ and $x>a, t>0$ continuous up to the boundaries and  
 satisfying the equation $iu_t+u_{xx}=0$.
 Now, integrating  (\ref{s}) by parts, it follows that  $u(x,t)=O\left(\frac{1}{x^2}\right)$ and $u_x(x,t)=O\left(\frac{1}{x}\right)$ as $|x|\to\infty$. 
Particularly,  $u(x,t)\in L^{1}(\mathbb{R})$ and $u_x(x,t)\to0$ as $|x|\to\infty$.

In order to prove that $u(x,t)$ satisfies the initial condition in (\ref{d2}), we notice that $u(x,0)$ reduces to 
$\frac{1}{2\pi}\int_{-\infty}^\infty e^{ikx}\hat u_0^{(1)}(k)dk$ for $x<a$
and to $\frac{1}{2\pi}\int_{-\infty}^\infty e^{ikx}\hat u_0^{(2)}(k)dk$ for $x>a$,
again applying Lemma \ref{r1} and Jordan's lemma for a circle in ${\mathbb C}_-$.

Finally,  the jump condition in (\ref{d2})
can be proven directly, starting from  (\ref{s}) and using the expressions (\ref{b}) for $u(a,t)$.

\section{The long-time asymptotics}
\label{ass}

The IST method, particularly, its realization as a Riemann--Hilbert problem method, has proven its high efficiency 
for studying asymptotic regimes  of nonlinear integrable equations, particularly, the long-time behavior of 
solutions of initial value problems \cite{DIZ,DZ} as well as initial boundary value problems (see, e.g., \cite{DP}),
where the result can be obtained by applying the so-called nonlinear steepest descent method for oscillatory Riemann--Hilbert problems.
For linear equations, as far as the solution is given in terms of contour integrals, it is the standard  steepest descent method that 
allows obtaining the asymptotic results. For problem (\ref{d2}), its application leads to the following

\begin{prop}
\label{asympt1}
Let $u_{0}(x)\in L^{1}({\mathbb R})$. Then the long-time asymptotics of the solution $u(x,t)$ 
of problem (\ref{d2}) along any ray $\frac{x}{2t}=\xi$ is as follows:
\begin{itemize}
\item if $q<0$, then, as $t\to\infty$, 
\begin{equation*}
u(x,t) = \frac{e^{it\xi^2-i\frac{\pi}{4}}}{2\sqrt{\pi t}}
\left(\hat{u}_0^{(1)}(\xi)-\frac{iqe^{-2i\xi a}\hat{u}_0^{(1)}(-\xi)-\xi\hat{u}_0^{(2)}(\xi)}{iq+\xi}\right)(1+o(1)),
\end{equation*}
for $x<a$ and
\begin{equation*}
u(x,t) = \frac{e^{it\xi^2-i\frac{\pi}{4}}}{2\sqrt{\pi t}}
\left(\hat{u}_0^{(2)}(\xi)-\frac{iqe^{-2i\xi a}\hat{u}_0^{(2)}(-\xi)+\xi\hat{u}_0^{(1)}(\xi)}{iq-\xi}\right)(1+o(1)),
\end{equation*}
for $x>a$.
\item if $q>0$, then as $t\to\infty$
\begin{equation*}
\begin{aligned}
u(x,t)&  =  qe^{q(x-a)+iq^{2}t}\left(\int_{-\infty}^{a}e^{q(y-a)}u_{0}(y)\,dy+
\int_{a}^{\infty}e^{-q(y-a)}u_{0}(y)\,dy\right)\\
&+\frac{e^{it\xi^2-i\frac{\pi}{4}}}{2\sqrt{\pi t}}
\left(\hat{u}_0^{(1)}(\xi)-\frac{iqe^{-2i\xi a}\hat{u}_0^{(1)}(-\xi)-\xi\hat{u}_0^{(2)}(\xi)}{iq+\xi}\right)(1+o(1)),
\end{aligned}
\end{equation*}
for $x<a$ and 
\begin{equation*}
\begin{aligned}
u(x,t)& =  qe^{-q(x-a)+iq^{2}t}\left(\int_{-\infty}^{a}e^{q(y-a)}u_{0}(y)\,dy+
\int_{a}^{\infty}e^{-q(y-a)}u_{0}(y)\,dy\right)\\
&+\frac{e^{it\xi^2-i\frac{\pi}{4}}}{2\sqrt{\pi t}}
\left(\hat{u}_0^{(2)}(\xi)-\frac{iqe^{-2i\xi a}\hat{u}_0^{(2)}(-\xi)+\xi\hat{u}_0^{(1)}(\xi)}{iq-\xi}\right)(1+o(1)),
\end{aligned}
\end{equation*}
for $x>a$.
\end{itemize}
\end{prop}

\section{Concluding remark}

In the present paper we have solved the IV problem for the linear Schr\"odinger equation with a point singular potential by the Fokas unified transform method. Although the initial and initial-boundary value problems for linear equations in $1+1$ dimension could be, in principle, treated by the classical methods  (particularly,
for the second-order equations), it is important to solve the original problem by
the unified transform method, since it is known as providing effective solutions to both linear and nonlinear integrable equations \cite{F, F2, P}. Therefore,
the present paper can be viewed as a step in attacking  (particularly, in a perturbative sense \cite{LF12, LF15})
the much more complicated initial value problem for the {\em nonlinear} Schr\"odinger equation with a point singular potential
$$
iu_{t}+u_{xx}+q\delta(x)u+2u|u|^{2}=0,
$$
in the case of general initial data $u_0(x)$ (not assuming it to be an even function). Moreover, the Fokas method is proved to be highly efficient for finding explicit solutions of some interface problems \cite{DS, DSS, DTV, DPS}. Therefore, in our paper we have established the applicability of the method for our particular interface problem.

\section{Acknowledgments}
The author would like to thank  D. Shepelsky for many useful conversations and helpful suggestions.
The partial support from the Akhiezer Foundation is gratefully acknowledged.
%

\end{document}